\documentclass{article}
\usepackage[utf8]{inputenc}

\usepackage{amsthm,amssymb,amsmath,array, authblk, forest,ltablex, subfigure, tabularx}
\usepackage{tikz}
\forestset{%
  gappy tree/.style={
    for tree={
      circle,
      draw,
      s sep'+=10pt,
      fit=band,
      align=center,
      calign angle=45,
      calign=fixed edge angles
    },
    },
  c phantom/.style={draw=none, no edge},
  c end/.style={draw=none, edge={dotted, shorten <=5pt}},
}
\usetikzlibrary{arrows,positioning, calc}
\tikzstyle{vertex}=[draw,fill=black!15,circle,minimum size=20pt,inner sep=0pt]

\begin{document}
\title{Hipster Tree Growth Rates}
\author{Alexander Wong}
\affil{University of Florida}
\date{\today}

\maketitle

  \subsubsection*{Abstract}

 A plane rooted tree is called a hipster tree if it has no nontrivial automorphisms. Equivalently, a tree is a hipster tree if no two siblings have isomorphic subtrees. We impose the hipster condition on various classes of rooted trees. By approximating the generating function for the number of such trees, we obtain bounds on their exponential growth rates. 
  \medskip

\section{Introduction}

Tree enumeration has been studied starting with Cayley's theorem on labelled free trees\cite{cayley}. Since then, both exact and asymptotic formulae for enumeration of many other classes of trees have been found\cite{moon}. In this paper we study asymptotic growth rates of classes of trees determined by constraining their automorphism group.

Automorphism groups of combinatorial objects are commonly studied in enumerative combinatorics as counting is often done up to symmetry. For trees, the distribution of sizes of the automorphism groups has been studied for various classes as in\cite{bona}\cite{olsson}. We instead consider trees with trivial automorphism groups (containing only the identity automorphism). We call such trees \textit{hipster trees}. Note that we may impose this hipster condition on any class of trees. We will focus on enumeration of subclasses of rooted plane trees with addition of this hipster condition. \\

We define a \textit{rooted plane tree} recursively as a tree with one special vertex $r$ called the root and the remaining vertices split into an ordered partition of rooted plane trees $\{T_1,T_2,\ldots,T_k\}$ where the root of each $T_i$ connected to $r$. We say $u$ is a \textit{parent} of $v$ if 
\begin{enumerate}
    \item There is an edge between $u$ and $v$.
    \item The path between $v$ and $r$ contains $u$.
\end{enumerate}
We also say $v$ is a \textit{child} of $u$. If two vertices share a parent, we say they are \textit{siblings}. \\

Let $T$ be a rooted plane tree. Let $v$ be a vertex of $T$. Let $S_v$ denote the subtree of $T$ consisting of all vertices weakly below $v$. Notice that a tree is hipster iff for any siblings $u,v$ we have $S_u \ncong S_v$. In this paper we will provide functional equations for the generating function of different classes of hipster trees. Then by approximating parts of the functional equation, we obtain generating functions that serve as bounds for the generating function of our class of hipster trees. Then by analyzing the singularities in our approximations, we obtain bounds on the exponential growth rate of the class of hipster trees.

\section{Binary Hipster Trees}
 We define a binary tree as a rooted plane tree where each vertex has at most two children and each child is distinguished as a left or right child even if it is the only child. Let $h_n$ be the number of binary hipster trees with $n$ vertices and let $H(x)$ be the generating function for $h_n$. Then $$H(x) = x((H(x))^2-H\left(x^2\right)+1)+1$$ as removing the root leaves two \textit{distinct} binary hipster trees. As evaluating generating functions at $x^2$ is difficult to handle, we can approximate $H(x)$ by approximating the $H\left(x^2\right)$ term. Let $\ell_n$ be the number of binary plane trees with $n$ vertices where each vertex has at most one child and let $L(x)$ be the associated generating function. Then $\ell_n = 2^{n-1}$ as each child can be either a left or right child. Then $$L(x) = \sum\limits_{i=0}^\infty 2^{i-1}x^i = 1+\frac{x}{1-2x}.$$
We have that $\ell_n \leq h_n$ for all $n$ as each tree counted by $\ell_n$ is a binary hipster tree. Let $G(x) = x(G(x)-L\left(x^2\right)+1)+1$ and let $g_n$ denote the coefficient of $x^n$ in the power series of $G(x)$. Then $G(x)$ is the generating function for binary trees where no two siblings induce the same subtree counted by $L(x)$. This is a weaker condition than the hipster condition as siblings can induce the same subtree as long as it is not counted by $L(x)$. Thus $g_n \geq h_n$. Solving for $G(x)$ gives us
\begin{align*}
    G(x) &= x(G(x)^2-L\left(x^2\right)+1)+1 \\
    0 &= xG(x)^2-G(x)-\frac{x^3}{1-2x^2}-x+x+1 \\
    G(x) &= \frac{1\pm \sqrt{1-4x(-\frac{x^3}{1-2x^2}+1)}}{2x}
\end{align*}
It is known that the exponential growth rate of the coefficients of a generating function is given by the reciprocal of the singularity closest to $0$\cite{flajolet}. The singularity closest to $0$ arises when $1+\frac{4x^4}{1-2x^2}-4x = 0$ which gives $$x = \frac{-2 - \sqrt{2} + \sqrt{14 + 4 \sqrt{2}}}{4}$$
Thus the exponential growth rate of $g_n$ is the inverse of the above singularity which gives $\approx 3.923909$. As $g_n \geq h_n$ this gives an upper bound for the exponential growth rate of $h_n$.\\\\
We now provide a lower bound the exponential growth rate of $h_n$. Let $c_n$ denote the number of binary trees on $n$ vertices and $C(x)$ the corresponding generating function. Then $c_n \geq h_n$. Note that the $c_n$ are the Catalan numbers and have generating function $\frac{1-\sqrt{1-4x}}{2x}$ \cite{stanley}. We now approximate the $H(x^2)$ term by $C(x^2)$ to get a lower bound for $h_n$. Let $$F(x) = x((F(x))^2-C\left(x^2\right))+1$$ and let $f_n$ be the coefficient of $x^n$ in $F(x)$. We claim $0 < f_n \leq h_n$. Proceed by induction. Note that $h_0 = f_0 = 1$. Assume $f_k\leq h_k$ for $k<n$. For any power series $G(x)$, let $[x^n](G(x))$ denote the coefficient of $x^n$ in $G(x)$. Consider $[x^n](F(x)) = [x^n](x((F(x))^2-C\left(x^2\right))+1)$ for $n>0$. We have
\begin{align*}
    [x^n](F(x)) &= [x^n](x((F(x))^2-C\left(x^2\right)+1)+1)\\
    f_n &= [x^{n-1}]((F(x))^2)-[x^{n-1}](C\left(x^2\right)+1) \\
    f_n &= \left(\sum\limits_{i=0}^{n-1} f_if_{n-i-1}\right)-c_{(n-1)/2}
\end{align*}
where we let $c_{(n-1)/2} = 0$ if $(n-1)/2$ is not an integer. Similarly we have 
\begin{align*}
    [x^n](H(x)) &= [x^n](x((H(x))^2-H\left(x^2\right)+1)+1)\\
    h_n &= [x^{n-1}]((H(x))^2)-[x^{n-1}](H\left(x^2\right)+1) \\
    h_n &= \left(\sum\limits_{i=0}^{n-1} h_ih_{n-i-1}\right)-h_{(n-1)/2}
\end{align*}
where we let $h_{(n-1)/2} = 0$ if $(n-1)/2$ is not an integer. Thus by the induction hypothesis and the fact that $c_n \geq h_n$ we have that $f_n \leq h_n$ as desired.\\
As $F(x) = x\left((F(x))^2-C\left(x^2\right)+1\right)+1$ we have
\begin{align*}
    F(x) &= x\left((F(x))^2-\frac{1-\sqrt{1-4x^2}}{2x^2}\right)+x+1 \\
    F(x)&=xF(x)^2-\frac{1-\sqrt{1-4x^2}}{2x}+x+1 \\
    0&= xF(x)^2-F(x)-\frac{1-\sqrt{1-4x^2}}{2x}+x+1 \\
    F(x) &= \frac{1\pm \sqrt{3-2\sqrt{1-4x^2}-4x^2-4x}}{2x}
\end{align*}
The singularity closest to $0$ arises when $3-2\sqrt{1-4x^2}-4x^2-4x = 0$ which gives $$x = \frac{1}{24} \left(-20 + (6400 - 768 \sqrt{69})^{1/3} + 4*2^{2/3} (25 + 3 \sqrt{69})^{1/3}\right).$$
Thus the exponential growth rate of $f_n$ is approximately the inverse of the above singularity which gives $\approx 3.923450$. As $h_n \geq f_n$ this gives a lower bound for the exponential growth rate of $h_n$. Thus the exponential growth rate of $h_n$ lies in $[3.923450,3.923909]$.

\section{Plane 1-2 Hipster Trees}
We can apply the same process for estimating the growth rates of other hipster trees provided the generating function of the total number of trees is sufficiently well understood. Recall a plane 1-2 tree is a rooted plane tree where each vertex has at most two children. Note that this is the same as a binary tree without the distinction between left and right children with no siblings. These trees are counted by the Motzkin numbers $m_n$ with generating function given by \cite{stanley}$$M(x) = \frac{1-x-\sqrt{1-2x-3x^2}}{2x^2}.$$
Once again, let $h_n$ be the number of plane 1-2 hipster trees with $n$ vertices and let $H(x)$ be the generating function for $h_n$. Then 
\begin{align*}
    H(x) &= x\left((H(x)-1)^2+(H(x)-1)-(H\left(x^2\right)-1)\right)+1\\
    &= x\left(H(x)^2-H(x)-H\right(x^2\left)+2\right)+1
\end{align*} 
as removing the root leaves two \textit{distinct} plane 1-2 hipster trees or a single plane 1-2 hipster tree (a single tree is double counted by the $H(x)^2$ term and subtracted once by the $-H(x)$ term). Let $\ell_n$ be the number of plane 1-2 trees with $n$ vertices where each vertex has at most one child and let $L(x)$ be the associated generating function. Then $$L(x) = \sum\limits_{i=0}^\infty x^i = \frac{1}{1-x}$$
as the only tree with $n$ vertices is a chain of length $n$. Then letting $G(x)$ be the generating function for plane 1-2 trees where no two siblings have the same tree in counted by $L(x)$ we have $$G(x) = x\left(G(x)^2-G(x)-L\left(x^2\right)+2\right)+1.$$ and $G(x)$ counts a superset of plane 1-2 hipster trees. Then
\begin{align*}
    G(x) &= x\left(G(x)^2-G(x)-L\left(x^2\right)+2\right)+1 \\
    G(x) &= x\left(G(x)^2-G(x)-\frac{1}{1-x^2}+2\right)+1 \\
    0&= xG(x)^2-(x+1)G(x)-\frac{1}{1-x^2}+2x+1
\end{align*}
and
$$G(x) = \frac{x+1\pm \sqrt{x^2+2x+1-4x(-\frac{1}{1-x^2}+2x+1)}}{2x}.$$
Numerical approximation gives a singularity at $x\approx 0.350277$ and the growth rate is $\approx 2.854882$. \\\\\
We now want to approximate $H\left(x^2\right)$ by $M\left(x^2\right)$. Let $$F(x) = x\left(F(x)^2-F(x)-M\left(x^2\right)+2\right)+1.$$ and consider $[x^n]F(x)$. We claim $0 < f_n \leq h_n$. Proceed by induction. Note that $h_0 = f_0 = 1$. Assume $f_k\leq h_k$ for $k<n$. We have 
\begin{align*}
    [x^n]F(x) &= [x^n]\left(x(F(x)^2-F(x)-M\left(x^2\right)+2\right)+1) \\
    f_n &= [x^{n-1}]\left(F(x)^2-F(x)-M\left(x^2\right)+2\right) \\
    f_n &= \left(\sum\limits_{i=0}^{n-1} f_if_{n-i-1}\right)-f_{n-1}-m_{(n-1)/2}+2\\
    f_n &= \left(\sum\limits_{i=1}^{n-1} f_if_{n-i-1}\right)-m_{(n-1)/2} +2 &&\textrm{as $f_0 = 1$}
\end{align*}
Similarly, considering $[x^n]H(x)$ we have
\begin{align*}
    [x^n]H(x) &= [x^n]\left(x(H(x)^2-H(x)-H\left(x^2\right)+2)+1\right) \\
    h_n &= [x^{n-1}]\left(H(x)^2-H(x)-H\left(x^2\right)\right) \\
    h_n &= \left(\sum\limits_{i=0}^{n-1} h_ih_{n-i-1}\right)-h_{n-1}-h_{(n-1)/2}+2\\
    h_n &= \left(\sum\limits_{i=1}^{n-1} h_ih_{n-i-1}\right)-h_{(n-1)/2} +2 &&\textrm{as $h_0 = 1$}
\end{align*}
By the induction hypothesis and the fact that $m_k\geq h_k$, we have that $f_n \leq h_n$ and the growth rate for $f_n$ gives us a lower bound for the growth rate for $h_n$. Then
\begin{align*}
    F(x) &= x\left(F(x)^2-F(x)-M\left(x^2\right)+2\right)+1 \\
    F(x) &= x\left(F(x)^2-F(x)-\frac{1-x^2-\sqrt{1-2x^2-3x^4}}{2x^4}+2\right)+1 \\
    0&= xF(x)^2-(x+1)F(x)-\frac{1-x^2-\sqrt{1-2x^2-3x^4}}{2x^3}+2x+1
\end{align*}
So 
$$F(x) = \frac{x+1\pm\sqrt{x^2+2x+1-\frac{2x^2-2+2\sqrt{1-2x^2-3x^4}}{x^2}-8x^2-4x}}{x^2}.$$
The singularity closest to $0$ arises when $$x^2+2x+1-\frac{2x^2-2+2\sqrt{1-2x^2-3x^4}}{x^2}-8x^2-4x=0.$$Numerical approximation gives $x\approx 0.354047$ giving a growth rate $\approx 2.824486$. Thus the exponential growth rate for $h_n$ lies in $[2.824486,2.854882]$.

\section{Binary Trees With Colored Right Edge}
We may choose other sequences that count trees and look at the subset of such trees that obey a hipster condition. Looking at the little Schr{\"o}der numbers, these count the number of binary trees where each right edge can be either blue or red \cite{stanley}. The little Schr{\"o}der numbers' generating function is given by $$S(x) = \frac{1-x-\sqrt{1-6x+x^2}}{4x}.$$
Once again, let $h_n$ be the number of binary hipster trees with right edges colored blue or red on $n$ vertices and let $H(x)$ be the generating function for $h_n$. Then $$H(x) = x\left((2H(x)^2- H(x) - 2H\left(x^2\right) + 2\right) + 1.$$ The $-H(x)$ term fixes the $2H(x)^2$ term which double counts trees where the root has no right child. The $+2,+1$ fix the double counting of the empty tree. 
We can bound $H(x)$. Let $\ell_n$ be the number of binary trees with right edge colored blue or red with $n$ vertices where each vertex has at most one child and let $L(x)$ be the associated generating function. We have that $\ell_n = 3^{n-1}$ as each edge can be left, blue and right, or red and right. Then $$L(x) = \sum\limits_{i=0}^\infty 3^{i-1}x^i = \frac{x}{1-3x}$$
as the only tree with $n$ children is a chain of length $n$. Then letting $$G(x) = x\left(2G(x)^2-G(x)-L\left(x^2\right)+2\right)+1$$ we have $G(x)$ is the generating function of a superset of the trees counted by $H(x)$ so $g_n \geq h_n$.
\begin{align*}
    G(x) &= x\left(2G(x)^2-G(x)-L\left(x^2\right)+2\right)+1 \\
    G(x) &= x\left(2G(x)^2-G(x)-\frac{x^2}{1-3x^2}+2\right)+1 \\
    0&= 2xG(x)^2-(x+1)G(x)-\frac{x^2}{1-3x^2}+2x+1
\end{align*}
and
$$G(x) = \frac{x+1\pm \sqrt{x^2+2x+1-4x(-\frac{x^2}{1-3x^2}+2x+1)}}{2x}.$$
Numerical approximation gives a singularity at $x\approx 0.174458$ and the growth rate is $\approx 5.732051$.\\\\
Let $$F(x) = x(2F(x)^2-F(x)-S\left(x^2\right)+2)+1.$$ Then
\begin{align*}
    [x^n]F(x) &= [x^n]\left(x(2F(x)^2-F(x)-S\left(x^2\right)+2)+1\right)+1) \\
    f_n &= [x^{n-1}]\left(2F(x)^2-F(x)-S\left(x^2\right)+2)+1\right) \\
    f_n &= \left(2\sum\limits_{i=0}^{n-1} f_if_{n-i-1}\right)-f_{n-1}-s_{(n-1)/2}+2\\
    f_n &= \left(2\sum\limits_{i=1}^{n-1} f_if_{n-i-1}\right)+f_{n-1}-s_{(n-1)/2}+2 &&\textrm{as $f_0 = 1$}
\end{align*}
Similarly, considering $[x^n]H(x)$ we have
\begin{align*}
    [x^n]H(x) &= [x^n]\left(x(2H(x)^2-H(x)-H\left(x^2\right)+2)+1\right)+1) \\
    h_n &= [x^{n-1}]\left(2H(x)^2-H(x)-H\left(x^2\right)+2)+1\right) \\
    h_n &= \left(2\sum\limits_{i=0}^{n-1} h_ih_{n-i-1}\right)-h_{n-1}-h_{(n-1)/2}+2\\
    h_n &= \left(2\sum\limits_{i=1}^{n-1} h_ih_{n-i-1}\right)+h_{n-1}-h_{(n-1)/2}+2 &&\textrm{as $h_0 = 1$}
\end{align*}
As $f_0=h_0$ and $h_i\leq s_i$, by induction we have that $f_i\leq h_i$ for all $i$.
\begin{align*}
    F(x) &= x(2F(x)^2-F(x)-S\left(x^2\right)+2)+1 \\
    F(x) &= x\left(2F(x)^2-F(x)-\frac{1-x^2-\sqrt{1-6x^2+x^4}}{4x^2}+2\right)+1 \\
    F(x) &= 2xF(x)^2-xF(x)-\frac{1-x^2-\sqrt{1-6x^2+x^4}}{4x^2}+2x+1\\
    0&= 2xF(x)^2-(x+1)F(x)-\frac{1-x^2-\sqrt{1-6x^2+x^4}}{4x^2}+2x+1
\end{align*}
So 
$$H(x) = \frac{x+1\pm\sqrt{x^2+2x+1-\frac{2x^2-2+2\sqrt{1-6x^2+x^4}}{x}-16x^2-8x}}{2x}.$$
The singularity closest to $0$ arises when $x^2+2x+1-\frac{2x^2-2+2\sqrt{1-6x^2+x^4}}{x}-16x^2-8x=0$. This happens at $x\approx 0.174465$ giving a growth rate of $\approx 5.731821$. Thus the exponential growth rate for $h_n$ lies in $[5.731821,5.732051]$.

\section{Further Directions}
There are many other classes of rooted plane trees this approach could be extended to (most notably the set of all rooted plane trees). A more difficult extension would be to some class of rooted labelled trees. Recall that a labelled rooted tree with $n$ vertices is a rooted tree with each of the numbers $1$ through $n$ associated with a vertex of the tree. We may define a labelled rooted hipster tree to be a labelled rooted tree where the underlying tree (without labels) is a rooted hipster tree. Consider the class of binary trees labelled such that each parent has a larger label than its children (called decreasing binary trees). Letting $D(x)$ be the \textit{exponential} generating function for such trees we get the differential equation $$D'(x) = D(x)^2-E(x)$$ where $E(x)$ is the exponential generating function for pairs of decreasing binary trees where the underlying trees are isomorphic and the labels are all distinct. Note that $E(x)$ itself is hard to describe, but perhaps it can be approximated in a convenient way. One may also ask whether the hipster condition could be defined differently for labelled trees such as insisting on subtrees of siblings having the same labels when considered as permutations.

\end{document}